\documentclass[10pt]{amsart}
\usepackage{amscd, amssymb}

\newcommand{\R}{{\mathbb R}}
\newcommand{\C}{{\mathbb C}}
\newcommand{\Z}{{\mathbb Z}}
\newcommand{\U}{{\mathcal U}}
\newcommand{\F}{{\mathcal F}}
\newcommand{\A}{{\mathcal A}}
\newcommand{\Prp}{{\mathcal P}}
\newcommand{\proj}{{\mathbb P}}
\newcommand{\E}{e}
\newcommand{\Of}{\Omega^{U, G}_*[\F]}
\newcommand{\Off}{\Omega^{U, G}_*[\F_1, \F_2]}
\newcommand{\wt}{\widetilde}
\newcommand{\Pn}{{\mathcal P}(\underline{n})}
\newcommand{\x}{\langle x \rangle}
\newcommand{\y}{\langle y \rangle}
\newcommand{\vphi}{\varphi}
\newcommand\mcC{{\mathcal{C}}}
\newcommand\mfp{\mathfrak{p}}
\newcommand\mcU{{\mathcal{U}}}
\newcommand\td{td}

\theoremstyle{plain}
\newtheorem{theorem}{Theorem}[section]
\newtheorem{proposition}[theorem]{Proposition}
\newtheorem{lemma}[theorem]{Lemma}
\newtheorem{corollary}[theorem]{Corollary}

\theoremstyle{definition}
\newtheorem{definition}[theorem]{Definition}
\newtheorem{example}{Example}

\theoremstyle{remark}
\newtheorem*{remark}{Remark}

\newcommand{\refT}[1]{Theorem~\ref{T:#1}}

\newcommand{\refP}[1]{Proposition~\ref{P:#1}}
\newcommand{\refD}[1]{Definition~\ref{D:#1}}
\newcommand{\refL}[1]{Lemma~\ref{L:#1}}
\newcommand{\refE}[1]{Equation~\ref{E:#1}}

\begin{document}
\title{The Geometry of the Local Cohomology Filtration in Equivariant Bordism}
\author{Dev P. Sinha}
\address{Department of Mathematics\\
          Brown University\\
          Providence, RI 02906}
\email{dps@math.brown.edu}
\subjclass{Primary: 55}

\begin{abstract}
We present geometric constructions which realize the local cohomology
filtration in the setting of equivariant bordism, with the aim of
understanding free $G$ actions on manifolds.  We begin by reviewing the
basic construction of the local cohomology filtration, starting with the
Conner-Floyd tom Dieck exact sequence.  We define this filtration
geometrically using the language of families of subgroups.  We then
review Atiyah-Segal-Wilson
$K$-theory invariants, which are well-suited for studying the manifolds
produced by our techniques.  We end by indicating potential applications
of these ideas.
\end{abstract}

\maketitle

\section{Introduction}
Local cohomology techniques are an area of rapid development in
equivariant topology.  Since their introduction by Greenlees
\cite{Gr1} they have been used for example to compute the $ku$-homology of
classifying spaces of abelian groups \cite{BrGr}, they have an
interesting relationship with the chromatic filtration \cite{GS}, and they
have been generalized beyond equivariant topology as announced by Dwyer
and Greenlees. They give a new duality between homology and cohomology of
classifying spaces, which can be played against the Universal Coefficient
duality.  See \cite{Gr2} for an excellent survey.

One purpose of this note is to better understand the local cohomology
filtration by casting it in the language of families.  We will see that
the local cohomology filtration is coarser than the
families filtration.  It is more manageable, however, because the algebra
which arises is that of localization.  Moreover, when one ``frees up''
the theory, it becomes more computable and gives rise to an interesting
relationship between the homology and cohomology, with coefficients in
any complex-oriented spectrum (for example ordinary  cohomology or
$K$-theory), of classifying spaces.

A strong appeal of equivariant bordism as a sub-field of equivariant
topology is in its connection with manifold theory.
Our geometric constructions allow one to prove collapse
results by constructing particular $G$-manifolds which represent classes
in the local cohomology spectral sequence.  Moreover, we may
also bring in Atiyah-Segal-Wilson localized
$K$-theory invariants to answer questions about differentials and
extensions.  On the geometric side, the constructions in this paper are in
part an answer to questions about realizing the  bordism of classifying
spaces geometrically.  These questions were posed to the author by
Botvinnik who has used such knowledge to solve cases of the
Gromov-Lawson-Rosenberg conjecture concerning spin manifolds which admit
positive scalar curvature metrics \cite{Bo, Stolz}.  Indeed, the main
ideas of this paper came out of conversations with Botvinnik and Sadofsky
while the author was visiting the University of Oregon.

This paper will blend well-established ideas of Conner-Floyd, tom Dieck,
Atiyah, Segal and Wilson with recent constructions of Greenlees and recent
insight of the author to show that the local cohomology filtration is
well-suited for addressing questions about both the algebra and geometry
associated to free $G$-actions on manifolds.  We begin the paper by giving
a retrospective of the Conner-Floyd-tom Dieck exact sequence in
$\Z/p$-equivariant bordism.  We then focus on the case of $(\Z/p)^2$,
defining the local cohomology filtration.   The new insights come when we
interpret this filtration  in the language of families in order to give a
geometric understanding of the filtration.  We then show how classes we
build using this filtration are well-suited to analysis through the use of
Atiyah-Segal-Wilson
$K$-theory invariants, taking our time to define these invariants. We
close with ways we see in which these techniques could be applied.

Note that Igor Kriz has told us that he has proved collapse of the local
cohomology spectral sequence for abelian groups using different methods.

\section{The Conner-Floyd-tom Dieck Exact Sequence}

Conner and Floyd were strongly interested in studying free $G$-actions
on manifolds \cite{CF}.  They saw that bordism theory, newly revolutionized
by Thom's work, would give information beyond what homology would see.
To that end, they were interested in bordism theory in which
representative classes and bordisms were equipped with a free $G$-action.
They realized that such a bordism theory was
equivalent to the bordism homology of the corresponding classifying
space.   Throughout we let $EG$ denote a contractible free $G$-space and
$BG$  denote the quotient of $EG$ by $G$.

\begin{proposition}\label{P:free}
The bordism module of stably complex free $G$-manifolds
is isomorphic to $\Omega^U_*(BG)$.
\end{proposition}

\begin{proof}
Consider the following diagram:
$$
         \begin{CD}
                 \wt{M}          @>\tilde{f}>>   EG     \\
                 @VVV                            @VVV    \\
                 M               @>f>>           BG.
         \end{CD}
$$
Given a representative $M$ with reference map $f$ to
$BG$, pull back the
principal $G$-bundle $EG$ to get $\wt{M}$, which
is in fact a free
$G$-manifold.  Conversely, starting with a free $G$-manifold
$\wt{M}$, there is no obstruction to constructing a map $\tilde{f}$
to $EG$.  Pass to quotients to obtain $f\colon M \to BG$.

These maps are well-defined, as we apply the previous argument to the
manifolds which act as bordisms.
The composites of these maps are clearly identity maps.
\end{proof}

Suffice it to say that the bordism module $\Omega^U_*(BG)$ contains a great
deal of information about free $G$ actions.  It is the main object of study
in this paper.  It remains mysterious for most groups, being easy to
describe only for cyclic groups \cite{Lan2} and having been computed at
all for only  a few groups, including elementary abelian groups
\cite{SWi}. Conner and Floyd realized that this module is best
understood in relation to bordism in which more general actions are
allowed.  We carefully define this theory now.


We begin
by recalling the definition of geometric complex equivariant bordism.  Fix
$\U_\C$ ($\U_\R$, respectively) a complex (resp. real) representation of
which a countably infinite direct sum of any
representation of $G$ appears as a summand.
Let $BU^G(n)$ ($BO^G(n)$, respectively) be the Grassmanian of complex (resp.
real) $n$-dimensional linear subspaces of $\U_\C$ (resp. $\U_\R$),
topologized as a direct limit of finite dimensional Grassmanians.

\begin{definition}
A tangentially complex $G$-manifold is a pair $(M, \tau)$ where $M$ is a
smooth $G$-manifold and $\tau$ is a lift to $BU^G(n)$ of the map to
$BO^G(2n)$ which classifies $TM \times \R^k$ for some $k$.
\end{definition}

We can define bordism equivalence in the usual way to get a geometric
version of equivariant bordism.

\begin{definition}
Let $\Omega^{U, G}_*$ denote the ring of tangentially complex
$G$-manifolds up to bordism equivalence.  Let $\Omega^{U,G}(X)$
denote the module (over $\Omega^{U,G}_*$) of tangentially complex
$G$-manifolds equipped with a reference map to $X$ up to bordism
equivalence.
\end{definition}

Recall that one may prove that geometric bordism theory has suspension
isomorphisms
or satisfies excision using transversality arguments \cite{MM}.
Because transversality holds between $G$-manifolds when the
range manifold has trivial action, $\Omega^{U,G}_*(-)$ has suspension
isomorphisms when the $G$-action on the suspension coordinate is trivial.
On the other hand, it does not have isomorphisms induced by
suspension by (one-point compactified) representations in general, so
we say that it is an equivariant homology theory ``indexed on the
trivial universe'' (see \cite{Wan} for the corresponding statement
in terms of spectra).

There is a map $\Omega^U_*(BG) \to \Omega^{U,G}_*$ defined simply because a
free $G$-manifold is a $G$ manifold.  Conner and Floyd realized that this
map fits in a natural exact sequence.  It will be best for us to describe
this exact sequence using the language of families, as we will be using this
language throughout the paper.  There is an excellent discussion of families
in \cite{Cos}.

\begin{definition}
A family of subgroups of a group $G$ is a set of subgroups $\F$
such that if $H \in \F$ and $H'$ is conjugate to a subgroup
of $H$ then $H' \in \F$.
\end{definition}

\begin{definition}
Let $\F$ be a family of subgroups of $G$.  An $\F$-space is a $G$-space
all of whose isotropy groups are in $\F$.
\end{definition}

There is a terminal object, unique up to homotopy,
in the category of $\F$-spaces which is called $E\F$.  In simpler terms,
$E\F$ is an $\F$-space which is contractible (after forgetting the
$G$-action).  For example, if $\F$ is the family consisting of only
the identity subgroup, $E\F$ is simply $EG$.
We will use the term $\F$-manifold to refer to a manifold which is an
$\F$-space.

If $\F_2 \subset \F_1$ are two families, define an $(\F_1, \F_2)$-manifold
to be $\F_1$-manifold with (possible empty) boundary such that its boundary
is an $\F_2$-manifold.

\begin{definition}
Two $(\F_1, \F_2)$-manifolds
$M$ and $N$ are bordant when there is an
$(\F_1)$-manifold $W$ such that $M \sqcup -N$ is equivariantly
diffeomorphic to a codimension zero sub-manifold of $\partial W$
and $\partial W - (M \sqcup N)$ is an $\F_2$-manifold.
\end{definition}

Let $\Of$ (respectively $\Off$)
denote the bordism module of tangentially complex $\F$ (resp.
$(\F_1, \F_2))$-manifolds.   If $\F_2 \subset \F_1$ are families, there is
an inclusion map $\Omega^{U, G}_*[\F_2] \to \Omega^{U, G}_*[\F_1]$ as well
as a boundary map $\Off \to \Omega^{U, G}_{*-1}[\F_2]$.  There are also
inclusion maps for pairs of families, and in particular
inclusion maps of the form
$\Omega^{U, G}_*[\F_1] \to \Omega^{U,G}_*[\F_1, \F_2]$, since
$\Omega^{U, G}_*[\F_1] \cong \Omega^{U, G}_*[\F_1, \phi]$, where
$\phi$ is the empty family.

\begin{proposition}\label{P:relfam}
Let $\F_3 \subset \F_2 \subset \F_1$ be three families of subgroups of $G$.
The sequence
$$ \cdots \to \Omega^{U, G}_*[\F_2, \F_3] \to \Omega^{U, G}_*[\F_1, \F_3] \to
\Omega^{U,G}[\F_1, \F_2] \overset{\partial}{\to} \Omega^{U, G}_{*-1}[\F_2,
\F_3] \to \cdots
$$ is exact.
\end{proposition}

This exact sequence is essentially the long exact sequence of a triple.
In fact, one may deduce it from the long exact sequence of a triple by
showing that $\Of \cong \Omega^{U,G}_*(E\F)$ and $\Off \cong \Omega^{U,
G}_*(E\F_1, E\F_2)$.

We now specialize to the case $G=\Z/p$, for which there are only
three families: the empty family, the family of the trivial (identity)
subgroup which we call $\E$, and the family of all subgroups which
we call $\A$.  We apply \refP{relfam} using these three families and
identify the terms.  Note that $\Omega^{U,G}_*[\A, \phi] = \Omega^{U,G}_*$.
Next consider $\Omega^{U,G}_*[\E, \phi]$.  By definition this is
the bordism module of free $G$-manifolds, which we have identified in
\refP{free}, and the map from it
to $\Omega^{U,G}_*$ is simply the forgetful map.
Hence the exact sequence reads as follows.

$$ \cdots \to \Omega^U_*(BG) \overset{i}{\to} \Omega^{U,G}_* \overset{j}{\to}
\Omega^{U,G}_*[\A,\E]
\overset{\partial}{\to} \Omega^U_{*-1}(BG) \to \cdots$$

Conner and Floyd went further in \cite{CF} to identify
$\Omega^{U,G}_*[\A,\E]$.  The map $j$ can be interpreted
in terms of ``reduction to fixed sets''.  The relevant general fact
is as follows.

\begin{proposition}\label{P:fixed}
An $(\F_1, \F_2)$-manifold $M$ is bordant to
any smooth neighborhood $\mathcal{N}(\bigcup_{H \notin \F_2} M^{H})$ of the
locus of points in $M$ fixed by
subgroups not in $\F_2$.
\end{proposition}

\begin{proof}
Let $W = M \times [0, 1]$, with ``straightened corners''.  Then $\partial
W$ is an $\F_2$ manifold outside of the codimension zero submanifolds $M
\times 0$ and a tubular neighborhood of $(\bigcup_{H \notin \F_2} M^{H})
\times 1$, so
$W$ is the required bordism.
\end{proof}

When $G = \Z/p$, tubular neighborhoods of fixed sets are diffeomorphic
to $G$ vector bundles over those fixed sets, so Conner and Floyd were
able to identify  $\Omega^{U,G}_*[\A,\E]$ in terms of bordism modules
of manifolds with trivial $G$-action which equipped with $G$-vector
bundles.

Returing to our original question, namely as to the structure of
$\Omega^U_*(BG)$, we see that this exact sequence gives rise to two-stage
filtration into the free $G$-manifolds which
are the boundaries of general $G$-manifolds (the cokernel of $\partial$)
and those which are not (the image of $i$).  We may now state that a main
goal of this paper is to give a similar description of the filtration
on $\Omega^U_*(BG)$ for more general $G$ which arises from Greenlees's
local cohomology machinery.
First, however, we wish to gain computational understanding from our
current exact sequence.

Unfortunately, to this day the geometric theories $\Omega^{U,G}_*$
are mysterious,
with partial knowledge for abelian $G$ and a complete
computation for $\Z/2$ which may be deduced from \cite{Si1}.  Hence
it is not computationally useful to study $\Omega^U_*(BG)$
using this exact sequence.  To make
this exact sequence more algebraically manageable requires two steps.

The first step was taken by tom Dieck in \cite{tD}, who realized that Conner
and Floyd's exact sequence was related to localization methods in equivariant
$K$-theory being developed at the time by Atiyah and Segal \cite{AS}.
To make this connection, tom Dieck defined a more homotopy theoretic
version of equivariant bordism.   He crafted a
spectrum $MU^G$, analogous to $MU$, whose corresponding infinite loop
space is
$\varinjlim_V \text{Maps} (S^{n\oplus V}, T(\xi^G_{|V|}))^G$, where $V$
ranges over isomorphism classes of complex representations of $G$,
$S^{n\oplus V}$ is the one-point  compactification of the Whitney sum of
$\C^n$ with trivial $G$ action and
$V$, and $T(\xi^G_{|V|})$ is the Thom space of the universal complex
$G$-bundle.   See \cite{Cos} for a recent treatment of tom Dieck's
construction.  Though this spectrum has been defined for thirty years and
studied actively, it is only recently that its coefficients have been
understood in any cases \cite{Kriz, Si1}.

The equivariant homology and cohomology theories corresponding to this
spectrum have suspension isomorphisms with respect to any representation.
Moreover, the cohomology theory has a Thom isomorphism for
complex $G$-vector bundles, and so has some theory of
characteristic classes.  An interesting class of $G$-vector bundles
are complex representations $V$, considered as bundles over a point.
The associated Euler classes are denoted $e_V \in MU_G^m(pt.) = MU^G_{-m}$.
There is a Pontryagin-Thom map from $\Omega^{U,G}_*$ to $MU^G_*$,
but it is not an isomorphism, since transversality arguments
fail in an equivariant setting \cite{tD} (in particular the Euler classes
$e_V$ are non-zero classes in negative homological degrees where
geometric bordism is zero by definition - see \cite{Cos}).

Reflecting on the Conner and Floyd exact sequence for $G= \Z/p$,  tom Dieck
considered the cofiber sequence
$$EG_+ \to S^0 \to \widetilde{EG},$$
whose long exact sequence in $\Omega^{U,G}$ theory is the Conner-Floyd
exact sequence.  We say that the resulting long exact sequence in $MU^G$
theory is the tom Dieck exact sequence.

tom Dieck realized that the map from $MU^G_*$ to $MU^G_*(\widetilde{EG})$
had an interpretation in terms of localization (and hence a connection
with the work of Atiyah and Segal \cite{AS}).  The key observation
is that one model for $EG$ for cyclic $G$ is as the unit sphere in
the representation $\oplus_\infty V$ which we denote $S(\oplus_\infty
V)$, where $V$ is the standard representation of $G$ as the roots of
unity in the complex numbers. Hence a model for $\widetilde{EG}$ is the
one-point compactification of that representation, which we denote
$S^{\infty V}$.

For any commutative ring $R$ and element $e \in R$ let $R[e^{-1}]$
denote the localization of $R$ obtained by inverting $e$.

\begin{lemma}\label{L:inveul}
For any $G$,
$\widetilde{MU^G_*}(S^{\oplus_\infty V}) \cong MU^G_*[e_V^{-1}]$ as rings.
\end{lemma}

\begin{proof}
The left-hand side $\widetilde{MU^G_*}(S^{\oplus_\infty V})$ is a ring
because $S^{\oplus_\infty V}$ is an $H$-space through the equivalence
$$S^{\oplus_\infty V} \wedge S^{\oplus_\infty V} \cong S^{\oplus_\infty V}.$$
To compute the left-hand side,
apply $\wt{MU^G_*}$ to the identification $S^{\oplus_\infty V} =
\varinjlim S^{\oplus_n V}$.  After applying the suspension
isomorphisms $\wt{MU^G_*}(S^{\oplus_k V}) \cong
\wt{MU^G_{*+|V|}}(S^{\oplus_{k+1}V})$, the maps in the resulting
directed system are multiplication by the $e_V$.
\end{proof}

tom Dieck also uses the fact that transversality arguments go through
in the presence of free $G$ actions to show that $MU^G_*(EG) \cong
MU_*(BG)$.  Here one could also use Adams's transfer argument
\cite{Ad2} to show $(MU^G \wedge EG_+)^G \cong MU \wedge BG_+$.  Hence the
tom Dieck exact sequence reads
$$
\cdots \to MU_*(BG) \to MU^G_* \to MU^G_*[{e_\rho}^{-1}] \to \cdots.
$$
Moreover, just as Conner and Floyd computed the third term in the
geometric  setting, one can compute the localization
$MU^G_*[{e_\rho}^{-1}]$. The fact that one can compute such a
localization is the starting point for the computation of
$MU^G_*$ for abelian $G$ \cite{Si1}, which  reveal that the structure of
these rings is quite complicated.  In fact, their structure is
complicated to an extent which renders the tom Dieck exact sequence
not usable for computation of $MU_*(BG)$.





The next step one must take to make this exact sequence more
computable is to ``free up'' the spectrum $MU^G$.  Consider
the mapping spectrum $\text{Maps}(EG_+, MU^G)$.  Its homotopy
groups are the homotopy groups of the $G$-maps from $EG_+$ to $MU^G$.
Because $EG_+$ is free, we have
$$\text{Maps}^G(EG_+, MU^G) = \text{Maps}(BG_+, i_*(MU^G))
= {\text{Maps}}(BG_{+} ,MU)$$
where $i_*$ takes a $G$-spectrum and passes to the underlying
spectrum.  Note as well that because $EG_+$ is simply $S^0$
non-equivariantly we have that $i_*({\text{Maps}}(EG_+, MU^G)) = MU$
as well.  Hence by Adams's transfer argument \cite{Ad2}
we have that
$$EG_+ \wedge_G {\text Maps}(EG_+, MU^G) = BG_+ \wedge i_*(\text
{Maps}(EG_+, MU^G)) = MU \wedge BG_+.$$
Finally, let $t^G MU^G$ denote the spectrum $\text{Maps}(EG_+, MU^G)
\wedge \widetilde{EG}$.  For cyclic $G$ the coefficients
of $t^G$ are obtained by inverting the Euler class of the standard
representation, as in \refL{inveul}, in the coefficients of $\text
{Maps}(EG_+, MU^G)$ which again are just $MU^*(BG)$.
Thus if we smash the spectrum $\text{Maps}(EG_+, MU^G)$ with
tom Dieck's cofiber sequence and take homotopy groups we get
$$
\cdots \to MU_*(BG) \to MU^*(BG) \to MU^*(BG)[{e_\rho}^{-1}] \to
\cdots
$$

The terms in this exact sequence are all computable.  Landweber used the
Gysin sequence for the fiber bundle $S^1 \to B\Z/p \to \C P^\infty$
to show in
\cite{Lan} that $MU^*(B\Z/p) \cong MU^*[[x]]/[p]_F(x)$ where
$[p]_F(x)$ is the $p$-series in the formal group law over $MU$ (see
\cite{Rav}).  More topologically, $p_F(x)$ is the Euler class of the $p$th
tensor power of the tautological bundle over $\C P^\infty$, whose
associated sphere bundle is the fiber bundle above.
Hence we may finally be explicit about the exact sequence, which reads
\begin{equation}\label{E:tateseq}
\cdots \to MU_*(BG) \to MU^*[[x]]/[p]_F(x) \to MU^*((x))/[p]_F(x) \to
\cdots,
\end{equation}
where in general we use $R((x))$ to denote $R[[x]][x^{-1}]$.
We may now compute $MU_*(BG)$, as was first done by Landweber \cite{Lan2}.
Note that $MU_* (BG)$ is a ring with unit, since it is isomorphic to
$MU^G_* (EG)$ and $EG$ is an $H$-space.

\begin{theorem}\label{T:zpcomp}
$MU_*(B\Z/p)$ is generated as a module by the unit class and by classes
$y_i$ with relations $p y_i + a_1 y_{i-1} + \cdots + a_{i-1} y_1 = 0$,
where $a_i$ is the coefficient of $x^i$ in the series $p_F(x)$.
\end{theorem}

\begin{proof}
We compute the cokernel and kernel of the map  $MU^*[[x]]/[p]_F(x) \to
MU^*((x))/[p]_F(x)$ in the sequence of Diagram~\ref{E:tateseq} above. The
cokernel is generated by the negative powers of $x$.  Let $y_i$ denote the
image of
$x^{-i}$ in $MU_*(BG)$.  By multiplying
the $p$-series by negative powers of $x$, we have that
$$p x^{-i} + a_1 x^{1-i} + \cdots + a_{i-1} x + {\text{non-negative
powers of x}} = 0,$$ so we deduce the stated relation.  Clearly all
relations arise from multiplying the $p$-series by a Laurent polynomial
in $x$, so these relations suffice to generate all relations.

The kernel is generated by the series $p_F(x)/x$.  It is a free
$MU_*$ module on one generator.

Because the cokernel consists entirely of $MU_*$-torsion ($p$-torsion, in
fact) and the kernel is free, we see that $MU_*(BG)$ splits as this sum
of these modules.  Moreover, since the ideal generated by the unit class
is non-zero and free, it must map isomorphically to its image under
$i$, namely the ideal generated by $p_F(x)/x$.
\end{proof}

It is interesting to note that the Atiyah-Hirzebruch spectral sequence
to compute this module has as many non-trivial extensions as possible.

To summarize the exact sequences we have introduced, for $G = \Z/p$ we
have a commutative diagram as follows.

$$
  \begin{CD}
   \Omega^U_*(BG) @>i>> \Omega^{U,G}_* @>j>> \Omega^{U,G}_*[\A,\E] \\
          @V{\cong}VV         @V{PT}VV                   @V{PT}VV    \\
      MU_*(BG) @>>> MU^G_* @>>> MU^G_*[{e_\rho}^{-1}]  \\
          @VVV         @VVV                   @VVV    \\
      MU_*(BG) @>>> MU^*(BG) @>>> MU^*(BG)[{e_\rho}^{-1}]\\
          @|        @|                  @|   \\
MU_*(BG) @>>> MU^*[[x]]/[p]_F(x) @>>> MU^*((x))/[p]_F(x)
  \end{CD}
$$

The middle two rows of this diagram are known as the Tate diagram
\cite{GM2}.

We may use this commutative diagram to give geometric representatives
of the classes in $\Omega^U_*(BG)$.  The key observation is the following.

\begin{proposition}[Lemma 3.1 of \cite{Cos}]
The class $e_\rho^{-1} \in MU^G_*[{e_\rho}^{-1}] $ is the image under
the Pontryagin-Thom map of $D(\rho)$, the unit disk in the representation
$\rho$, which is a $G$-manifold with $\A$-boundary.
\end{proposition}

\begin{proof}[Sketch of Proof]
The appropriate setting for the proof is that of ``stable manifolds''.
Given a manifold with boundary $M, \partial M$ and a map of pairs from
this manifold to $D(V), S(V)$ for some representation $V$ there is a
Pontryagin-Thom map which produces an element of
$\pi_{\text{dim} M - \text{dim} V}(MU^G)$.  Non-equivariantly, such a map
does not produce any classes beyond those produced by the standard
Pontryagin-Thom map (it cannot, since the usual Pontryagin-Thom map is
an isomorphism).  Equivariantly, this procedure does produce new classes,
in fact giving all of $MU^G_*$ \cite{BH}.

For example, $e_V$ is the image under this map of the stable manifold in
which $M$ is a point which maps to $0 \in D(V)$.  It is then conceivable
that the inverse of this class in the localization is given by the
image under the Pontryagin-Thom map of $D(V)$ mapping to a point
(zero-disk).
\end{proof}

\begin{corollary}
The generator $y_i \in \Omega^U_{2i-1}$ is represented by
$S(\oplus_i \rho)$.
\end{corollary}

It is quite straightforward to geometrically understand the unit class in
$MU_*(BG)$ as well by running through the isomorphism of
Proposition \ref{P:free}.

\begin{proposition}\label{P:repunit}
The unit class in $MU_*(BG)$ is represented by $G$ itself, as a
zero-dimensional manifold.
\end{proposition}

Hence, we not only have an algebraic computation of $\Omega^U_*(BG)$ but
we have explicit geometric understanding of most of the generators as
(quotients of) free $G$-manifolds which are boundaries of arbitrary
$G$-manifolds. Note that having a free $G$-manifold as the boundary of an
arbitrary
$G$-manifold is precisely the setting for studying the free $G$-manifold
using the $K$-theoretic invariants of Atiyah and Segal \cite{AS} as
further developed by Wilson \cite{Wi}.  We will develop these ideas below.

\bigskip

Our goal for more general $G$ is to develop similar understanding of
$\Omega^U_*(BG)$, with as much control as possible over both the algebra
and the geometry.  We will outline a program for doing so in this paper,
describing a computational tool which generalizes the discussion of this
section, namely the local cohomology filtration. There are in fact two
ways in which one may generalize the exact sequences developed in this
section, the other known as the families filtration, which is described
very well in \cite{Cos}.  Briefly, one filters $EG$ by a sequence of maps
$EG \to E\F_1 \to \cdots E\F_n$, where $F_{i+1}$ is obtained from $F_i$
by adding one (conjugacy class of) subgroup.  For many theories such as
bordism, the equivariant homology of the cofiber of $E\F_i \to
E\F_{i+1}$ is equivalent to (non-equivariant) homology of some sort of
classifying space (much as in Proposition \ref{P:free}).  The  local
cohomology filtration is coarser than this filtration but surprisingly
much more computable once one ``frees up'' the theories involved.

\section{The Local Cohomology Filtration}

Just as the survey \cite{Car} successfully illustrated equivariant
theory by concentrating on $\Z/2$, we are going to focus on $(\Z/p)^2$,
giving an indication of how things would generalize as we go along.

In the last section it was useful for compuations to notice that the third
term Conner-Floyd-tom Dieck sequence, namely the one associated to
$S^{\oplus_\infty \rho}$, is a localization.  To continue to have
localizations involved in our computations, we are going to assemble
$EG$ from pieces of the form $S^{\oplus_\infty V}$ for various $V$.

For $G = (\Z/p)^2$, $EG$ is no longer given by an action on a sphere
but on a product of spheres.  Let $G$ be generated by elements $x$ and $y$
(so that $px = py = x + y - x - y = 0$).  Let $V_x$ (respectively $V_y$)
be the standard representation of $\Z/p$ pulled back to $G$ through its
projection to $G/\langle x \rangle$ (respectively $G/\langle y
\rangle$).   We have that
$$EG = S(\oplus_\infty V_x) \times S(\oplus_\infty V_y).$$
Hence we smash the cofiber sequence $S(\oplus_\infty V_x)_+ \to S^0
\to S^{\oplus_\infty V_x}$ with the correspondence sequence with $V_x$
replaced by $V_y$ to get a commutative diagram as follows:

\begin{equation}\label{E:sqr}
  \begin{CD}
    EG_+ @>>> S(\oplus_\infty V_x)_+ @>>>  S(\oplus_\infty V_x)_+ \wedge
                                             S^{ \oplus_\infty V_y}\\
          @VVV         @VVV                   @VVV    \\
    S(\oplus_\infty V_y)_+  @>>> S^0 @>>> S^{ \oplus_\infty V_y} \\
          @VVV         @VVV                   @VVV    \\
     S(\oplus_\infty V_y)_+ \wedge S^{ \oplus_\infty V_x}  @>>>
         S^{ \oplus_\infty V_x}   @>>> S^{ \oplus_\infty (V_x \oplus V_y)}.
  \end{CD}
\end{equation}

In this situation, $\Sigma^2 EG_+$ is known as the total cofiber of the
square diagram consisting of the last two spaces of the last two rows in
the above diagram.  All of the spaces in question, except for $S^0$, are
of the form $S^{\oplus_\infty V}$ as desired.

We pause for a moment to recall
basic definitions of  cubical diagrams and their total cofibers (from
\cite{Go}) and give a canonical filtration on those total cofibers.



\medskip

Let $\underline{n} = \{1, \cdots, n\}$ and let $\Pn$ denote the category
of subsets of $\underline{n}$  with morphisms given by inclusion.  A
cubical diagram of spaces (of dimension $n$) is a functor from $\Pn$ to
based spaces.  We often let $X_\bullet$ denote such a functor, $X_S$
denote the value of this functor on $S \in \text{ob}(\Pn)$ and $f_{S,
S'}$ denote the unique morphism from $X_S$ to $X_{S'}$ when $S \subseteq
S'$.  A morphism of  cubical diagrams is defined as usual for a diagram
(or functor) category.

Cubical diagrams are a convenient diagram category in part because a
morphism of cubical diagrams may be viewed as a cubical diagram of
dimension greater by one.  Let $i$ be an element of $\underline{n+1}$
and by let $\vphi_i$ denote the order-preserving inclusion from
$\underline{n}$ to $\underline{n+1}$ for which $i$ is not in the image.
We may define an $n+1$-dimensional cubical diagram
$Z_\bullet$ from a morphism $X_\bullet \to Y_\bullet$ by letting
$Z_{\vphi_i(S)} = X_{S}$, $Z_{\vphi_i(S) \cup i} = Y_S$ and the map from
$Z_{\vphi_i(S)}$ to $Z_{\vphi_i(S) \cup i}$ be the map from $X_S$ to
$Y_S$.  Conversely, given a cubical diagram of dimension $n$ one may view
it as a morphism of cubical diagrams of dimenstion $n-1$ in $n$ ways.

We now define the total cofiber of a cubical diagram, which is a functor
from cubical diagrams to spaces we will use extensively.  For a cube
$I^n$ whose vertices are naturally labelled by subsets of $\underline{n}$
let $\partial_S I^{n}$ for a subset $S$ of $\underline{n}$ denote the
face whose vertices are subsets $S$.

\begin{definition}
The total cofiber of an $n$-dimensional cubical diagram $X_\bullet$ is
the quotient of the union $\coprod X_S \times \partial_{\underline{n} -
S} I^n$  through the identifications
\begin{itemize}
\item For $x \in X_S$ and $S \subset S'$, $x \times
\partial_{\underline{n} - S'} I^{n} \subset x \times
\partial_{\underline{n} - S} I^{n}$ is identified with $f_{S, S'}(x)
\times  \partial_{\underline{n} - S'} I^{n}$.
\item Points of the form $* \times \partial_A I^n$, where $*$ is the
basepoint of $X_S$ and $A$ is a face which does not meet the initial
vertex (labelled by the empty set), are all identified to a basepoint.
\end{itemize}
\end{definition}
 

There are Puppe sequences for cubical diagrams.  Given a map of cubes
$X_\bullet \to Y_\bullet$ one may take the Puppe sequence for each
map $X_S \to Y_S$.  These sequences fit together to define a sequence of
cubical diagrams which we also call a Puppe sequence.  Moreover, we
have the following.

\begin{proposition}\label{P:puppe}
Let $C_k$ be the total cofiber of the $k$th cube in the Puppe
sequence of a map $X_\bullet
\to Y_\bullet$.  The sequence $C_0 \to C_1 \to \cdots$ is a sequence
of spaces weakly equivalent to the Puppe sequence associated to the
induced map from the total cofiber of $X_\bullet$ to the total cofiber of
$Y_\bullet$.
\end{proposition}

There is a canonical filtration on the total cofiber of a cubical
diagram.

\begin{proposition}\label{P:ss}
Let $h$ be a homology theory and let $X_\bullet$ be a cubical diagram of
dimension $n$.  There is an $n$-column spectral sequence whose $E_1$ term
is given by $$E^1_{p,q} = \bigoplus_{\substack{S \subset
\underline{n}\\ \#S = n-p}} h_{q-\#S}(X_S)$$ and
$$d_1|_{h_*(X_S)} = \Sigma_{S' \supset S}
    (-1)^{\#S'\cdot \mu(S,S')} (f_{S,S'})_*,$$
where $\mu(S, S')$ is the element of $S'$ not in $S$,
which converges to the homology of the total cofiber of $X_\bullet$.
\end{proposition}

\begin{proof}[Partial proof]
We supply a proof, which is immediately generalizable, in the case of
square diagrams.

Consider the sequence of square diagrams
$$
\begin{CD}
* @>>> * \\
@VVV       @VVV\\
* @>>> X_{1,2}
\end{CD}
\overset{f_1}{\to}
\begin{CD}
* @>>> X_1 \\
@VVV       @VVV\\
X_2 @>>> X_{1,2}
\end{CD}
\overset{f_2}{\to}
\begin{CD}
X_\phi @>>> X_1 \\
@VVV       @VVV\\
X_2 @>>> X_{1,2}.
\end{CD}
$$
These maps give rise to a sequence of maps between total cofibers of
those diagrams, which we denote $T_0 \overset{f_1}{\to} T_1
\overset{f_2}{\to} T_2$.  As there is for any such sequence of maps, there
is a three-column spectral sequence whose
$E_1$ term consists of the homology of $T_0$, the cofiber of $f_1$ which
we denote $cof(f_1)$, and $cof(f_2)$.  Let $g_i$ denote the map from
$T_i$ to $cof(f_i)$ (with $g_0$ the identity map by convention) and
$\partial_i$ the map in the Puppe sequence  from $cof(f_i)$ to $\Sigma
T_{i-1}$.  The $d_1$ of this spectral sequence is given by $g_{i-1} \circ
\partial_i$.  This spectral sequence converges to the homology of
$T_2$.

We have that $T_0$ is simply $X_{1,2}$.  Next we see that the cofiber of
the map from
$T_0$ to $T_1$ is the total cofiber of the cube defined by $f_1$ as a
map of squares.  This total cofiber is in turn the total cofiber of the
square
$$
\begin{CD}
{\text{cofiber}} (* \overset{f_1}{\to} *) @>>>
                      {\text{cofiber}} (* \overset{f_1}{\to} X_1) \\
@VVV       @VVV\\
{\text{cofiber}} (* \overset{f_1}{\to} X_2) @>>>
                 {\text{cofiber}} (X_{1,2} \overset{f_1}{\to} X_{1,2}),
\end{CD}
=
\begin{CD}
* @>>> X_1\\
@VVV     @VVV\\
X_2 @>>> *
\end{CD}
$$
which is simply $\Sigma X_1 \wedge \Sigma X_2$.  Similarly, the cofiber of
the map from $T_1$ to $T_2$ is $\Sigma^2 X_\phi$.

Finally, we identify the $d_1$ differential.  As we stated above,
it suffices to understand the boundary maps in the Puppe sequence
of total cofibers.  By \refP{puppe}, we may instead look at the
Puppe sequence of cubes.  Unraveling the definitions, we see that the
boundary maps used to define $d_1$ in a spectral sequence of a filtration
are a wedge of suspensions of structure maps from the original cube, with
signs introduced by the interchange of suspension coordinates.

\end{proof}

We apply this proposition to the square diagram
$$
  \begin{CD}
    S^0 @>>> S^{ \oplus_\infty V_y} \\
       @VVV                   @VVV   \\
     S^{ \oplus_\infty V_x}   @>>> S^{ \oplus_\infty (V_x \oplus V_y)}.
  \end{CD}
$$
where the homology theory in question is either $MU^G$ or respectively
$\text{Maps}(EG_+, MU^G)$, and make identifications using \refL{inveul} in
order to deduce the following theorem originally due to Greenlees. Once
again we restrict our attention to the case of
$(\Z/p)^2$.

\begin{theorem}[Greenlees]\label{T:locss}
There are three-column spectral sequences converging to
$MU_*(B({\Z/p}^2))$ whose $E^1$ terms are given by
$$  MU^G_*[e_{V_x}^{-1}, e_{V_y}^{-1}] \overset{+,-}{\longleftarrow}
   MU^G_*[e_{V_x}^{-1}] \oplus MU^G_*[e_{V_y}^{-1}]
     \overset{}{\longleftarrow} MU^G_*,$$
and
\begin{multline}\label{E:locss} MU_*[[x,y]][(xy)^{-1}]/ (p_F(x), p_F(y))\\
\overset{+,-}{\longleftarrow}
    MU_*[[x,y]][x^{-1}]  / (p_F(x), p_F(y)) \oplus
          MU_*[[x,y]][y^{-1}]/ (p_F(x), p_F(y))\\
     \overset{}{\longleftarrow} MU^*(BG) = MU_*[[x,y]]/ (p_F(x)),
p_F(y)),
\end{multline}
where the grading defined so that these columns are in degrees $-2$,
$-1$ and $0$ and where the maps above are the canonical maps arising from
the fact that each of the columns, which are graded rings, is obtained
from the next by localization.
\end{theorem}

This spectral sequence is called the local cohomology spectral sequence,
because in fact the $E^2$ term is precisely what is known as the local
cohomology of the ring $MU^G_*$, respectively $MU^*(BG)$, at the ideal
generated by the Euler classes $e_{V_x}$ and $e_{V_y}$ \cite{Har}.  The
algebra of local cohomology and its relevance to this situation is
well-documented in the work of Greenlees and his collaborators \cite{Gr1,
Gr2, GM}.  In short, local cohomology groups of a ring $R$ at an ideal
$I$ are the derived functors of the $I$-torsion functor.   Because of its
connection with well-developed homological algebra, the local cohomology
filtration has lead to many advances in equivariant topology, and the
topology of classifying spaces in particular.

\section{A Geometric Version of the Local Cohomology Filtration}

Here we will take a geometric
approach by recasting this spectral sequence in the language of families.
First we replace our diagram of nine spaces from Diagram \ref{E:sqr} -
the cofiber sequence of cofiber sequences - with a diagram of bordism
modules defined by families.  Then, we give geometric definitions of the
differentials in this spectral sequence.  Such
analysis leads to conditions for collapse of the spectral sequence, which
we verify for $G = (\Z/p)^2$.

Let $\E$ and $\A$ continue to denote the families (of subgroups of
$(\Z/p)^2$ now) consisting of only the trivial group and of all subgroups,
respectively.  By abuse, let $\x$ denote the family consisting of the
subgroup generated by $x$ along with the trivial subgroup, and similarly
for $\y$.  Let $\x \cup \y$ denote the family consisting of the
subgroups generated by $x$ and by $y$ along with the trivial subgroup.

Consider the following diagram of bordism modules.

\begin{equation}\label{E:famsquare}
  \begin{CD}
\Omega^{U,G}_*[\E, \phi] @>h_{11}>> \Omega^{U,G}_*[\x, \phi]
	@>h_{12} >>	\Omega^{U,G}_*[\x, \E] \\
          @Vv_{11}VV         @Vv_{12}VV                   @Vv_{13}VV  \\
    \Omega^{U,G}_*[\y, \phi] @>h_{21}>> \Omega^{U,G}_*[\A, \phi]
		@>h_{22}>> \Omega^{U,G}_*[\A, \y] \\
          @Vv_{21}VV         @Vv_{22}VV               @Vv_{23}VV \\
     \Omega^{U,G}_*[\y, \E] @>h_{31}>> \Omega^{U,G}_*[\A, \x] @>h_{32}>> 
	\Omega^{U,G}_*[\A, \x \cup \y].
  \end{CD}
\end{equation}

There are boundary maps which continue the sequences defined by the first
two rows and columns.  For example, $h_{13} \colon
\Omega^{U,G}_*[\x,\E] \to
\Omega^{U,G}_{*-1}[\E, \phi]$ by defined by taking the boundary of
a representative class.  These are simply the boundary maps of
\refP{relfam}. There are also boundary maps for the last column and row.
Given an
$(\A, \x \cup \y)$-manifold $M$ representing a class in bordism, the
boundary map of the last column sends $[M]$ to $[\nu((\partial M)^x)]$ the
class represented by a tubular neighborhood (in $\partial M$) of the
$x$-fixed set of $\partial M$.

\begin{proposition}
The rows and columns of Diagram \ref{E:famsquare}, continued with
the boundary maps defined above, form long exact sequences.
\end{proposition}

For the first and second rows and columns this is an application
\refP{relfam}.  For the last row and column, the verification is similar
to that of the long exact sequence of a pair in bordism.

As mentioned above, we are constructing an analog of Diagram \ref{E:sqr}.

\begin{proposition}
The sequences of Diagram \ref{E:famsquare} map to those of
Diagram~\ref{E:sqr} through the Pontryagin-Thom map.
\end{proposition}

Our cubical diagram formalism allows us to use Diagram~\ref{E:famsquare}
to construct an analog of the local cohomology spectral sequence.

\begin{theorem}\label{T:geomss}
There is a spectral sequence whose $E^1$ term is given by
$$ \Omega^{U,G}_*[\A, \x \cup \y] \overset{v_{23}-h_{32}}{\longleftarrow}
   \Omega^{U,G}_*[\A, \y] \oplus \Omega^{U,G}_*[\A, \x]
     \overset{h_{22} \oplus v_{22}}{\longleftarrow} \Omega^{U,G}_*[\A,
\phi] = \Omega^{U,G}_*,$$ which converges to $MU_*(BG)$.  This spectral
sequence maps to those of \refT{locss}.
\end{theorem}

Because the maps in this spectral sequence are geometrically defined we
may further our understanding of how the various sub-quotients in the
local cohomology filtration give rise to classes in
$MU_*(BG)$.

For example, we
see that the module $\Omega^{U,G}_*[\A,\x \cup \y]$ maps to
$\Omega^{U,G}_*[\E,\phi] = MU_*(BG)$ by sending a class $[M]$ to
$[\partial(\nu(\partial M)^x)]$, which is the boundary of a tubular
neighborhood (in $\partial M$) the the $x$-fixed set of $\partial M$.
We may verify this by starting at $(\A, \x \cup \y)$ and composing a
verticle arrow with a horizontal arrow in Diagram \ref{E:famsquare}.  If
$M_1 = \partial W_1$ and $M_2 = \partial W_2$ are free $\Z/p$-manifolds
which bound arbitrary $\Z/2$-manifolds, then $W_1 \times W_2$ is
naturally a $(\Z/p)^2$-manifold with boundary which represents a class in
$\Omega^{U,G}_*[\A, \x \cup \y]$.  Under the map to $MU_*(BG)$, this
class maps to $[M_1 \times M_2]$.

We may understand the $d_2$ differential as the map defined by $h_{12}$
applied to a lift of a class $[M] \in \Omega^{U,G}_*[\A, \phi]$ which
maps to zero under $v_{22}$ and $h_{22}$.  Geometrically, this
differential determines whether a manifold $[M]$ which is cobordant to
both a manifold with no
$\x$-fixed points and to another with no $\y$-fixed points is actually
cobordant to a free manifold.


Of course, in order to prove that a differential is zero it suffices to
find non-zero classes in $MU_*(BG)$ to which the classes in the spectral
sequence coincide.

\begin{theorem}
The spectral sequence of Diagram~\ref{E:locss} in \refT{locss} collapses
at $E^2$.
\end{theorem}

\begin{proof}
The kernel of the $d_1$ map
$$ MU_*[[x,y]]/ (p_F(x)) \to MU_*[[x,y]][x^{-1}]  / (p_F(x), p_F(y))
                \oplus           MU_*[[x,y]][y^{-1}]/ (p_F(x), p_F(y))
$$ is the ideal generated by the series $p_F(x)/x \cdot p_F(y)/y$.  We
claim that this series is in the image of the map from $MU_*(BG)$ to this
quotient in the filtration, and hence cannot support a
$d_2$ differential.  For dimensional reasons, no other differentials are
possible.

The class which maps to the series $p_F(x)/x \cdot p_F(y)/y$ is the unit
class in $MU_*(BG)$, represented by $G$ itself, a zero-dimensional
manifold.  This follows from the fact that the unit class in
$MU_*(B\Z/p)$ mapped to $p_F(x)/x$ in $MU^*(B\Z/p)$ in the Tate sequence,
as we established in \refT{zpcomp}.
\end{proof}

Perhaps the most fun and interesting aspect of this analysis occurs for
in the middle column of the spectral sequence, even though it cannot
support a differential simply for dimensional reasons.
If two classes $[M] \in \Omega^{U,G}_*[\A, \y]$ and $[N] \in
\Omega^{U,G}_*[\A, \x]$ have the same image in $\Omega^{U,G}_*[\A,\x \cup
\y]$ then if we take $[\partial M] \in \Omega^{U,G}_*[\y, \phi]$ it must
go to zero under $v_{21}$ and hence lift to $\Omega^{U,G}_*[e, \phi] =
MU_*(BG)$.
 
Geometrically, becuase $[M]$ and $[N]$ have the same image in
$\Omega^{U,G}_*[\A,\x
\cup \y]$ it means that there is a cobordism $W$ between a tubular
neighborhood of $\bigcup_{H \neq \x {\text{or}} \y} M^H$, which we
denote $\nu_M$, and a tubular neighborhood of $\bigcup_{H \neq \x
{\text{or}} \y} N^H$, which we call $\nu_N$.
The free $G$-manifold to which $\partial M$ is cobordant may be
constructed by glueing together $\partial M - \nu((\partial M)^y)$,
$\partial N - \nu((\partial N)^x)$ and $S(\nu(\partial W)^x)$.

It is easy to be more explicit if we suppose now that the tubular
neighborhoods $\nu_M$ and $\nu_N$ are diffeomorphic and that we can extend
that diffeomorphism to include a neighborhood of $(\partial M)^y$ in $M$
and $(\partial N)^x$ in
$N$.  Then one may plumb $M$ and $N$ together by identifying these
neghborhoods (and ``rounding the corners'' in a canonical way) to define a
manifold with boundary whose boundary has a free $G$-action.  This
boundary represents the coset in $MU_*(BG)$ associated to $[M]$ and
$[N]$ in this filtration as above.

For example, let $G = {\Z/2}^2$.
Given a complex representation $W$ let $\proj(W)$ denote the
space of complex one-dimensional subspaces of $W$ with inherited
$G$-action.  We consider $\proj(\C \oplus V_x)$, where by abuse
$\C$ denotes the one-dimensional trivial representation, which is just
the Riemann sphere where $x$ acts trivially and $y$ acts by
multiplication by $-1$.  Then $\proj(\C \oplus V_x) \times D(V_y)$
represents a class in $\Omega^{U,G}_4[\A, \y]$ and, similarly,
$\proj(\C \oplus V_y) \times D(V_x)$
represents a class in $\Omega^{U,G}_4[\A, \x]$.  These classes map to the
same class in $\Omega^{U,G}_4[\A,\x \cup \y]$ since neighborhoods of
their $(x+y)$-fixed sets are both diffeomorphic to two copies of
$D(V_x \oplus V_y)$.  These diffeomorphisms which do extend to the
boundaries as needed above.  If we plumb $\proj(\C \oplus V_x) \times
D(V_y)$ and $\proj(\C \oplus V_y) \times D(V_x)$ along these
neighborhoods we get a four-manifold $P$ whose boundary is free.  It is
useful to picture the real analog in which two bands - $S^1 \times D^1$ -
are glued along $S^0
\times I^2$ to get a surface diffeomorphic to $S^2$ with four open disks
removed.

We will show in the next section the the class $[\partial P] \in
MU_3(B(\Z/2)^2)$ is non-zero.  In fact, it must then repsent a ``Tor
class'', by which we mean a class which is not in the image
of the tensor product of $MU_*(B\Z/2)$ with itself under the K\"unneth
map.  As far as we know, ours is the first construction of a
representative for a class in $MU_*(BA)$ for any abelian group $A$ which
is not a union of products of spheres.  Our
plumbing constructions greatly increase the tools with which one can
create free $A$-actions.

\medskip

Our techniques clearly generalize beyond ${\Z/p}^2$.  One may chase
through the analog of Diagram~\ref{E:famsquare} to define classes generated by
processes composed the basic ones we have used in this section: taking
boundaries, taking tubular neighborhoods of fixed sets, and using
cobordisms between equivalent classes.

\section{Atiyah-Segal-Wilson Invariants}

In this section we describe some powerful invariants which may be used
to study free $G$-manifolds which are boundaries of arbitrary
$G$-manifolds, as are most manifolds constructed from the geometric local
cohomology filtration.  These invariants are essentially characteristic
numbers in localized equivariant $K$-theory.  We will take a geometric
approach, though it will be clear how to stabilize as suggested to us by
May.

Recall that in ordinary non-equivariant homology, characteristic numbers
encode the homomorphism
$$ \chi \colon MU_*(X) \to H_*(BU \times X),$$
which sends the fundamental class in homology of $M$ to its image under
$(\tau \times f)_*$ in $H_*(BU \times X)$.  When $X$ is an $H$-space (in
particular, when $X$ is a point), this homomorphism is a ring
homomorphism.  There is a corresponding construction at the level of
spectra (see for example \cite{Stong}).

\begin{proposition}\label{P:chiboard}
The homomorphism $\chi$ corresponds with the Boardman homomorphism
$$h \colon \pi_*(MU \wedge X_+) \to \pi_*(MU \wedge H \wedge X_+).$$
\end{proposition}

\begin{proof}
Let $\nu$ be the normal bundle of $M$ embedded in some sphere $S^{k+N}$, and
consider the following commutative diagram
\begin{equation}\label{E:board}
    \begin{CD}
	H_{k+N}(S^{k+N}) @>c_*>> H_{k+N} T(\nu)   @>T(\tau \times f)_*>>
		H_{k+N}(T(\xi_N) \wedge X_+)  \\
	&&   @A{\text{Thom} \cong}AA	@A{\text{Thom} \cong}AA  \\
	&&		H_k(M)	@>(\tau \times f)_*>>	H_k(BU(N) \times X),
    \end{CD}
\end{equation}
where $c$ is the collapse map onto $\nu$ as in the Pontrijagin-Thom
construction.  Then the top composite is the homology Boardman
homomorphism, and the bottom map is $\chi$.  Commutativity of the right
square and the fact that $c_*$ is a homology isomorphism in dimension
$k+N$ gives the desired correspondence between homomorphisms.
\end{proof}

\begin{remark}
When using the natural basis for $\pi_*(MU \wedge H
\wedge X_+)$ these numbers correspond to the characteristic numbers of
the stable normal bundle, not tangent bundle, of $M$.
\end{remark}

A more homotopy theoretic point of view of characteristic numbers
facilitates their definition for generalized, possibly equivariant,
cohomology theories which have a Thom isomorphism for complex
vector bundles.  Given $M$ a stably complex manifold, $E$ a cohomology
theory for which the Thom isomorphism theorem holds,
and $x \in E^*(M)$ we can define ``evaluation of $x$ on the
fundamental class of $M$'', as motivated by \refE{board}.

\begin{definition}
Let $x \in E^*(M)$.  With notation as in \refP{chiboard},
define the number $x[M]$ to be the image of $x$ under the composite
$$
\wt{E}^*(M_+) \overset{\cong}{\to}
	\wt{E}^{*+|\nu|}(T(\nu)) \overset{c^*}{\to}
		\wt{E}^{*+|\nu|}(S^V) \overset{\cong}{\to} E_{*-n}(pt.),
$$
where $n$ is the dimension of $M$.
\end{definition}

This composite map goes by many names in the literature, including
``Gysin map'' and ``wrong-way map''. The class $x[M]$ is a bordism
invariant of $M$ when
$x$ is a characteristic class of the tangent bundle of $M$.  As in the
case of ordinary cohomology, these numbers encode a ring homomorphism.
To avoid discussion of the fundamental class and facilitate computation,
we may define this ring homomorphism by giving explicit formulae for the
characteristic classes dual to multiplicative generators of $E_*(BU)$.
We define this homomorphism now in the case of equivariant $K$-theory,
which for finite complexes may be defined in an entirely analogous way to
ordinary $K$-theory by using $G$-vector bundles (see
\cite{Segal} for a full discussion).

\begin{remark}
Equivariant $K$-theory characteristic numbers play an important role in
geometry, by work of Atiyah, Segal and Singer \cite{AS}.  For example
if $\gamma$ is a holomorphic vector bundle over a complex $G$-manifold $M$,
let $H^i(M, \mathcal{O}(\gamma))$ denote the complex $G$-vector space which
is the cohomology of the sheaf of holomorphic sections of $\gamma$.  Then
the equivariant Riemann-Roch theorem of \cite{AS} states that
$$ \sum_i (-1)^i H^i(M, \mathcal{O}(\gamma)) = \gamma[M].$$
This sheaf cohomology will often vanish for $i>0$ so we will have that
$\gamma[M]$ is isomorphic to the $G$-vector space of holomorphic sections
of $\gamma$.
\end{remark}

\begin{definition}
Let $\gamma$ be an $n$-dimensional complex $G$-vector bundle over $X$, and let
$\gamma$ also denote the corresponding class in $K_G^0(X)$.
The characteristic class $\beta_i(\gamma)$ is the coefficient of $t^i$ in the
power series
$$  \log \left( \sum_{i \geq 0} \lambda^i(\gamma - n) \cdot t^i \right)
\in \wt{K}_G^0(X)[[t]].$$
\end{definition}

Recall the coefficient ring for equivariant $K$-theory, which satisfies a
strong form of Bott periodicity.
$K_G^0$ is isomorphic to the the representation ring $R(G)$, while
$K_G^1$ is zero.

\begin{definition}\label{D:kappa}
Let $M$ represent a class in $\Omega^{U,G}_n$.  Define the ring homomorphism
$\kappa \colon \Omega^{U, G}_* \to R(G)[[b_i]]$ by setting
$$ \kappa([M]) = \sum_{I = (n_1, n_2, \cdots n_k)} \beta_I(TM)[M]\cdot b_I$$
where $\beta_I(E) =  \prod_{1 \leq i \leq k} \left( \beta_i(E) \right)^{n_i}$
and $b_I = b_1^{n_1} b_2^{n_2} \cdots b_k^{n_k}$.
\end{definition}

The homomorphism $\kappa$ extends to the
Boardman homomorphism $\pi_*(MU^G) \to \pi_*(MU^G \wedge K^G)$.

We think of the coefficients of $\kappa$ as invariants of $\Omega^{U, G}_*$
with values in the representation ring.  A standard approach to the
representation ring is through character theory.  We can ask what parts
of the structure of a $G$-manifold $M$ we need to retain in order to
understand the value of $E[M]$ at a single conjugacy class of $G$.
For example, the value of $E[M]$ at the identity
element, namely its dimension,
will only depend on the underlying manifold $M$ and not the $G$-action.
Atiyah and Segal show in \cite{AS} that the characters associated to these
invariants  at a single conjugacy class depend only on the fixed-set
structure of $M$ at that conjugacy class, which is given by the image
of $[M]$ in $\Omega^{U,H}_*[\A, \Prp]$, where $H$ is the cyclic subgroup
of $G$ generated by a representative of the conjugacy class and $\Prp$ is
the family of proper subgroups of $H$.

Atiyah and Segal noticed that the geometry of reduction to fixed sets
coincided in algebra to localization, which we have already seen in
\refL{inveul}.

\begin{definition}
Given a conjugacy class $\mcC$ of the representation ring $R(G)$ let
$\mfp(\mcC)$ be the prime ideal of elements of $R(G)$ whose characters do
not vanish at $\mcC$.
Let $R(G)_{\mfp(\mcC)}$ denote the localization of
$R(G)$ at $\mfp(\mcC)$, and let $\lambda \colon R(G) \to
R(G)_{\mfp(\mcC)}$ denote the canonical map of $R(G)$ into this
localization.
\end{definition}

\begin{definition}
Let $ev_{\mcC} \colon R(G)_{\mfp(\mcC)} \to \C$
denote the map which gives the value of a localized character at $\mcC$.
This map is  well-defined because characters of denominators in
$R(G)_{\mfp(\mcC)}$ do not  vanish on $\mcC$.
\end{definition}

\begin{theorem}[Atiyah-Segal]\label{T:aslocal}
Let $\mcC$ be a conjugacy class of $G$ and let $h \in \mcC$.  Then the
composite
$$
\Omega^{U, G}_* \overset{\kappa}{\to} R(G)[[b_i]] \overset{\lambda}{\to}
   	R(G)_{\mfp(\mcC)}[[b_i]] \overset{ev_{\mcC}}{\to} \C[[b_i]]
$$
factors though the homomorphism $\Omega^{U, G}_* \to
\Omega^{U,H}_*[\A, \Prp]$. In other words, there is a homomorphism
$$\kappa_{\mcC} \colon  \Omega^{U,H}_*[\A, \Prp] \to \C[[b_i]]$$ such
that
$\kappa_{\mcC} \circ \vphi^h = ev \circ \lambda \circ \kappa$, where
$\vphi^h$ is the canonical map from $\Omega^{U,G}_*$ to
$\Omega^{U,H}_*[\A, \Prp]$.
\end{theorem}

We will give a formula for $\kappa_{\mcC}$ below.

Now note that the maps from $\Omega^{U,G}_*$ to $\Omega^{U,H}_*[\A, \Prp]$
all factor through $\Omega^{U,G}_*[\A, \E]$.  This observation in
conjunction with \refT{aslocal} leads to powerful invariants for
$MU_*(BG)$.  Consider the following diagram
\begin{equation}\label{E:wilson}
  \begin{CD}
   \Omega^{U,G}_*[\A, \phi] @>\alpha >> \Omega^{U,G}_*[\A, \E]
				@>\beta>> 
	\Omega^{U,G}_{*-1}[\E, \phi] = MU_{*-1}(BG)\\
          @V{\kappa}VV         @V{\kappa'}VV           \\
     R(G)[[b_i]]  @>\chi>> \bigoplus_\mcC \C[[b_i]],
  \end{CD}
\end{equation}
where by abuse we are now using $\chi$ to denote the character map in
representation theory. By exactness of the first row and commutativity of
the square we may deduce the following.

\begin{theorem}\label{T:kinvars}
If class $[M, \partial M] \in \Omega^{U,G}_*[\A, \E]$ has image under
$\kappa'$ which is not in the image of $\chi$ then $[\partial M]$ is
non-zero $MU_{*-1}(BG)$.
\end{theorem}

In \cite{Wi}, G. Wilson shows that the converse of this theorem is true
for groups with periodic cohomology.

We now present the cohomological definition of $\kappa_{\mcC}$, as given
in \cite{AS}.  First note that by \refP{fixed} an element $[M]$ of
$\Omega^{U,H}_*[\A, \Prp]$ is represented by a tubular neighborhood of
$M^H$, which by the equivariant tubular neighborhood theorem is
diffeomorphic to the total space of an $H$ vector bundle over $M^H$.  We
will be viewing classes as represented by such bundles for the purposes of
defining our fomulae.   We first define an equivariant version of the
Chern character for $G$-bundles over trivial $G$-spaces.

\begin{definition}
Let $E$ be a complex $G$-vector bundle over a trivial $G$-space $X$, which
decomposes as $\oplus_{\rho_i \in Irr^(G)} E_{\rho_i}$ (all such
$G$-bundles do so - see \cite{Segal}). For each $g \in G$ define
$$ch(E)(g) = \sum_{\rho_i \in Irr(G)} \chi_{\rho_i}(g) \cdot ch(E_{\rho_i}),$$
where $\chi_{\rho_i}(g)$ is the trace of $g$ acting on $\rho_i$.
\end{definition}

We next introduce the Todd class which is needed to translate from
$K$-theory to cohomology.   Let $c_i$ denote the $i$th
Chern class of a complex vector bundle $\gamma$ in
$H^{2i}(X)$.  Form the ring
$$ H^*(X)[a_1, \cdots, a_n]/(s_i(a_1, \cdots, a_n) = c_i),$$
where $s_i$ denotes the $i$th elementary symmetric polynomial in the
variables indicated.  Thus, any symmetric polynomial in the $a_i$ will
be a polynomial in the $c_i$.

\begin{definition}\label{D:todd}
Let $M$ be a unitary manifold.  The Todd class of $M$, denoted $\td(M)$,
is defined by taking $\gamma$ to be a complex vector bundle which is a
lift of the stable tangent bundle of $M$.
Then set
$$ \td(M) = \prod_{i=1}^{n} \frac{a_i}{1-e^{-a_i}}.$$
\end{definition}

Now let $H$ be a cyclic group of order $n$, and fix a generator $h$ of $H$.
Let $\zeta = e^{\left(\frac{2 \pi i}{n} \right)}$.
And let $\rho_i$ be the representation of $H$ in which the generator
$h$ acts on $\C$ by multiplication by $\zeta^i$.

\begin{definition}
Let $E_{\rho_i}$ be a complex $H$-vector bundle over a base with trivial
$H$-action and fiber isomorphic to $\oplus_k \rho_i$ for some $k$.  Define
$$
\mcU_i(E_{\rho_i}) = \prod_{1 \leq j \leq k} \left( \frac{\zeta^i - 1}
					{\zeta^i - e^{-a_j}} \right).
$$
If $E = \oplus_{1 \leq j < n} E_{\rho_i}$.
Define
$$
\mcU(E) = \prod_i \mcU_i(E_{\rho_i}).
$$
\end{definition}

We are now ready to present the cohomological formula we will use
to define the maps $\kappa_{\mcC}$.

\begin{definition}\label{D:loceval}
Let $H$ be a cyclic group, generated by $h$.
Let $X \in \Omega^{U,H}_*[\A,\Prp]$ be represented by a $G$-bundle whose
base is
$M$, a stably complex manifold, whose total space we denote $\nu$ and
whose fiber is isomorphic to a representation $V$ such that $V^H = 0$.
Let $E$ be an element of $K_H^0(M)$. Define
$$
   E[X] = \left\{ \frac{ch(E)(h) \cdot \mcU(\nu) \cdot \td(M)}
		{\det(1 - h | V^*)} \right\} [M]
$$
\end{definition}

\begin{definition}
Let $\mcC$ be a conjugacy class of $G$.  Let $h \in \mcC$, and let $H$ be
the cyclic subgroup generated by $h$.
Let $X$ represent a class in $\Omega^{U,H}_n[\A, \Prp]$.  Define the ring
homomorphism
$\kappa_{\mcC} \colon \Omega^{U,H}_*[\A, \Prp] \to
R(G)_{\mfp(\mcC)}[[b_i]]$ by
$$ \kappa_{\mcC}([X]) =
    \sum_{I = (n_1, n_2, \cdots n_k)} \beta_I(TX-n)[X] \cdot b_I$$
where $\beta_I(E) =  \prod_{1 \leq i \leq k} \left( \beta_i(E) \right)^{n_i}$
and $b_I = b_1^{n_1} b_2^{n_2} \cdots b_k^{n_k}$.
\end{definition}

Though the formula in \refD{loceval} looks complicated, it is manageable
in some cases, especially when the fixed sets are isolated points.

\begin{example}\label{E:compX1}
Let $G = \Z/p$ and fix a generator $g$.  Consider the class
$X_{1, \rho_1} \in \Omega^{U,G}_2[\A, \E]$ represented by $D(\rho_1)$.
In the notation of
\refD{loceval}, for $X_{1, \rho_1}$ we have
$V = \rho_1$, $\nu$ is a trivial bundle, and $M$ is a point.
Substituting, $\mcU(\nu) = \td(M) = 1$, so that
$$\kappa_{g}(X_{1, \rho_1}) = \frac{1}{1 - \zeta^{-1}} \cdot f(b_i),$$
where $f$ is a power series in $b_i$ whose coefficients lie in $\Z[\zeta]$
and whose constant term is one.
\end{example}

\begin{example}
Let $G = \Z/4$ and fix a generator $g$.  Consider the complex
$G$-manifold  with free boundary $X$
whose underlying topology is that of two complex disks, where the generator $g$
acts by mapping the disks to one another so that $g^2$ acts by multiplying
by $-1$ on each disk.  Then since $g$ and $g^3$ do not fix any points of
$X$, we have $\kappa_g(X) = \kappa_{g^3}(X) = 0$.  Restricted to
$\langle g^2 \rangle$, $X$ is two
copies of $X_{1, \rho_1}$, where $\rho_1$ is the sign representation of
$\Z/2$.  By the above example $\kappa_{g^2}(X)$ is integral with
constant term one.  Thus, from the constant terms of the $\kappa_{\mcC}(X)$
we construct a localized character whose values are $(0, 1, 0)$ at
$g$, $g^2$ and $g^3$, respectively.  Even though each entry of this
localized character is an integer, the character itself is not integral in
that it is not in the span of $(1, 1, 1)$, $(i, -1, -i)$, $(-1, 1, -1)$
and $(-i, -1, i)$, the values of characters of irreducible 
representations of $G$
at $g$, $g^2$ and $g^3$.  Thus $X$ is not the fixed set of any complex
$\Z/4$ surface.
\end{example}

Finally, we apply these invariants to the ${\Z/2}^2$ manifold $P$ defined
at the end of the previous section.  $P$ has two fixed points, whose
normal bundles are $V_x \oplus V_y$ in the notation of the previous
section.  If we look at the conjugacy class of $x+y$, there are only
these two fixed points and the computations proceed as in the first
example above to give a value of $1/4 + 1/4 = 1/2$ multiplied by some
integral power series in $b_i$, which shows that the boundary of $P$
represents a non-trivial class in $MU_3(B{\Z/2}^2)$.  At the conjugacy
class $x$ (respectively
$y$) the fixed set is a $\proj^1$, whose Todd genus is $1 +
c_1(\proj^1)/2 = 1+a$, where
$a$ is the generator of $H^2$.  As in the first example, there is a two
in the denominator coming from the action of $x$ on the normal bundle,
and so we once again get that the value of the localized character is
$1/2$ multiplied by an integral series.

\section{Directions for Further Work}

We leave with a couple suggestions for projects which may use these
techniques.

\begin{itemize}

\item Do computations in the local cohomology spectral sequence for
$MU_*(BA)$.  Compare this filtration to the K\"unneth filtration (they
seem to be very closely related), which is one of the main tools for
computations so far
\cite{SWi}.  Atiyah-Segal-Wilson invariants may be useful in solving
extension problems.

\item Translate these constructions to the Spin setting in order to
approach the Gromov-Lawson-Rosenberg conjecture for finite groups
\cite{Bo, Stolz}.  Part of the approach of \cite{Bo} is to find arbitrary
$G$-manifolds whose boundaries generate free $G$-bordism.

\end{itemize}

\end{document}